\documentclass[11pt]{article}
\usepackage{amssymb,amsfonts,amsmath,latexsym,epsf,tikz,url}

\newtheorem{theorem}{Theorem}[section]
\newtheorem{proposition}[theorem]{Proposition}

\newtheorem{conjecture}[theorem]{Conjecture}
\newtheorem{corollary}[theorem]{Corollary}
\newtheorem{lemma}[theorem]{Lemma}

\newcommand{\proof}{\noindent{\bf Proof. }}
\newcommand{\qed}{\hfill $\square$\medskip}

\begin{document}

\title{Symmetry breaking in planar and maximal outerplanar graphs }

\author{
Saeid Alikhani $^{}$\footnote{Corresponding author}
\and
Samaneh Soltani 
}

\date{\today}

\maketitle

\begin{center}
Department of Mathematics, Yazd University, 89195-741, Yazd, Iran\\
{\tt alikhani@yazd.ac.ir, s.soltani1979@gmail.com}
\end{center}


\begin{abstract}
The distinguishing number (index) $D(G)$ ($D'(G)$) of a graph $G$ is the least integer $d$
such that $G$ has a vertex (edge) labeling   with $d$ labels  that is preserved only by a trivial automorphism.  In this paper we consider the  maximal outerplanar  graphs (MOP graphs) and show that MOP graphs, except $K_3$,   can be distinguished by
at most two vertex (edge) labels.   We also  compute the distinguishing number and the distinguishing index of Halin and   Mycielskian graphs.
 
\end{abstract}


\medskip
\noindent{\bf AMS Subj.\ Class.}: 05C25 

\section{Introduction}
Let $G = (V,E)$ be a simple, connected and undirected graph, and let ${\rm Aut}(G)$ be its   automorphism group.
     A labeling of $G$, $\phi : V \rightarrow \{1, 2, \ldots , r\}$, is said to be \textit{$r$-distinguishing},  if no non-trivial  automorphism of $G$ preserves all of the vertex labels. The point of the labels on the vertices is to destroy the symmetries of the
graph, that is, to make the automorphism group of the labeled graph trivial.
Formally, $\phi$ is $r$-distinguishing if for every non-identity $\sigma \in {\rm Aut}(G)$, there exists $x$ in $V$ such that $\phi(x) \neq \phi(\sigma(x))$. The \textit{distinguishing number} of a graph $G$ is defined  by
\begin{equation*}
D(G) = {\rm min}\{r \vert ~ G ~\text{{\rm has a labeling that is $r$-distinguishing}}\}.
\end{equation*} 

This number has defined in \cite{Albert}.  The \textit{distinguishing index} $D'(G)$ of a graph $G$ is the least number $d$ such that $G$ has an edge labeling  with $d$ labels that is preserved only by the identity automorphism of $G$.  The distinguishing edge labeling was first defined by  Kalinowski and Pil\'sniak \cite{R. Kalinowski and M. Pilsniak} for graphs. Obviously, this invariant is not
defined for graphs having $K_2$ as a connected component. If a graph has no nontrivial automorphisms, its distinguishing number is  $1$. In other words, $D(G) = 1$ for the asymmetric graphs.
 The other extreme, $D(G) = \vert V(G) \vert$, occurs if and only if $G = K_n$. The distinguishing index of some examples of graphs was exhibited. For 
 instance, $D(P_n) = D'(P_n)=2$ for every $n\geqslant 3$, and 
 $D(C_n) = D'(C_n)=3$ for $n =3,4,5$,  $D(C_n) = D'(C_n)=2$ for $n \geqslant 6$. It is easy to see that the value $|D(G)-D'(G)|$ can be large. For example $D'(K_{p,p})=2$ and $D(K_{p,p})=p+1$, for $p\geq 4$.

 A \textit{maximal outerplanar graph} is an outerplanar graph that cannot have any additional edges added to it while preserving outerplanarity. Every maximal outerplanar graph with $n$ vertices has exactly $2n - 3$ edges, and every bounded face of a maximal outerplanar graph is a triangle. Every maximal outerplanar graph satisfies a stronger condition than Hamiltonicity: it is \textit{pancyclic}, meaning that for every vertex
 $v$ and every $k$ in the range from three to the number of vertices in the graph, there is a length $k$ cycle containing $v$. A cycle of this
 length may be found by repeatedly removing a triangle that is connected to the rest of the graph by a single edge, such that the
 removed vertex is not $v$, until the outer face of the remaining graph has length $k$.
 In \cite{arvind2}, Arvind et al. designed
efficient algorithms for computing the distinguishing numbers of trees and planar graphs.  Arvind et al. proved that the distinguishing number of a planar graph
can be computed in time polynomial in the size of the graph, \cite{arvind}.
 Fijav{\v{z}} et al.  showed that
every 3-connected planar graph is 5-distinguishing colorable except $K_{2,2,2}$ and $C_6 +\overline{K_2}$, \cite{Fijav}. They also proved that every 3-connected bipartite planar graph is 3-distinguishing colorable except $Q_3$
and $R(Q_3)$.

\begin{figure}[ht]
 	\begin{center}
 		\includegraphics[width=0.4\textwidth]{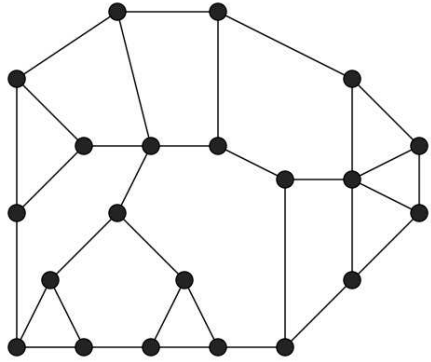}
 		\caption{\label{hal-hal3}Examples of a Halin  graph.}
 	\end{center}
 \end{figure}

 A \textit{Halin graph} is a type of planar graph, constructed by connecting the leaves of a tree into a cycle. The tree must have at least four vertices, none of which has exactly two neighbors; it should be drawn in the plane so none of its edges cross (this is called planar embedding), and the cycle connects the leaves in their clockwise ordering in this embedding. Thus, the cycle forms the outer face of the Halin graph, with the tree inside it, see Figure \ref{hal-hal3}. Every Halin graph is a Hamiltonian graph, and every edge of the graph belongs to a Hamiltonian cycle. Moreover, any Halin graph remains Hamiltonian after deletion of any vertex \cite{cornuejols1983halin}. 
 
\begin{figure}
	\begin{center}
		\includegraphics[width=0.8\textwidth]{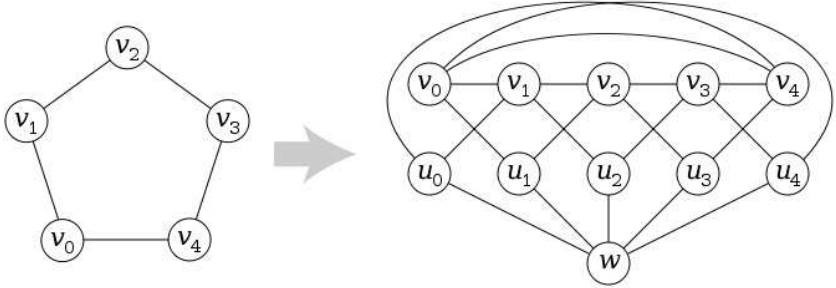}
		\caption{\label{gorzqu}The Gr\"otzsch graph as the Mycielskian of a 5-cycle graph.}
	\end{center}
\end{figure}

  The \textit{Mycielskian} or \textit{Mycielski graph}  of an undirected graph is a larger graph formed from it by a construction of Jan Mycielski (1955) \cite{Mycielski}. The construction preserves the property of being triangle-free but increases the chromatic number; by applying the construction repeatedly to a triangle-free starting graph, Mycielski showed that there exist triangle-free graphs with arbitrary large chromatic number.
    Let the $n$ vertices of the given graph $G$ be $v_0$, $v_1,..., v_{n-1}$. The Mycielski graph $\mu(G)$ of $G$ contains $G$ itself as an isomorphic subgraph, together with $n+1$ additional vertices: a vertex $u_i$ corresponding to each vertex $v_i$ of $G$, and another vertex $w$. Each vertex $u_i$ is connected by an edge to $w$, so that these vertices form a subgraph in the form of a star $K_{1,n}$. In addition, for each edge $v_iv_j$ of $G$, the Mycielski graph includes two edges, $u_iv_j$ and $v_iu_j$. The illustration shows Mycielski's construction as applied to a 5-vertex cycle graph with vertices $v_i$ for $0 \leq i \leq 4$. The resulting Mycielskian is the \textit{Gr\"otzsch graph}, an 11-vertex graph with 20 edges, see Figure \ref{gorzqu}. Thus, if $G$ has $n$ vertices and $m$ edges, $\mu(G)$ has $2n+1$ vertices and $3m+n$ edges.

\medskip

We study  the  distinguishing number and the distinguishing index of planar and maximal outerplanar graphs in Section 2. Then we study symmetry breaking for the  Halin graphs which are planar graphs, in Section 3. Finally, in Section 4, we obtain  the  distinguishing number and the distinguishing index of  Mycielski graphs.

\section{Maximal outerplanar graph}

In this section we study symmetry breaking for planar and maximal outerplanar graphs. We first  consider planar graphs.  

\subsection{Planar graphs} 

\begin{theorem}\label{ubboundplana}
	If  $G$ is  a connected graph with clique number four  and the maximum degree $\Delta \geq 5$, then $D(G)\leq \Delta -1$.
\end{theorem}
\proof
Let $v_0,v_1,v_2,v_3$ be vertices of an induced subgraph $K_4$ in $G$. We label the vertex $v_0$ with label $\Delta -1$, and vertices $v_1$, $v_2$ and $v_3$ with labels 1, 2 and 3, respectively. Next we label the remaining unlabeled vertices in $N_G(v_0)$ with different labels in the set $\{1, 2, \ldots , \Delta -2\}\setminus \{3\}$. Now we consider the BFS tree $T$ of $G$ with root $v_0$, hence $N_G(v_0)=N_T(v_0)$. We label the vertex $v_0$ and each vertex in $N_T(v_0)$, exactly the same as $G$. For labeling of the remaining vertices of $T$, we use the following mathematical induction.\\
Let $x$ be  a labeled vertex of $T$ at distance $i$ from $v_0$ where $i\geq 1$. Then, we label the adjacent vertices to $x$ which are at distance $i+1$ of $v_0$, with different labels in the set $\{1, 2, \ldots , \Delta -1\}$.

We now consider this labeling for vertices of $G$. If there is an  induced    subgraph $K_4$ in $G$ with vertices of labels exactly the same labels as the induced subgraph $v_0v_1v_2v_3$, i.e., labels 1,2,3 and $\Delta -1$, except  the induced subgraph $v_0v_1v_2v_3$, then we change the labels of vertices  of this induced subgraph by permuting of labels of vertices in $N_T(x)$ where $x$ is some of vertices of that induced subgraph $K_4$. So without loss of generality we can assume that the induced subgraph $K_4$ with vertices $v_0,v_1,v_2,v_3$ is the only induced subgraph $K_4$ in $G$ with label set $\{1,2,3, \Delta -1\}$. Hence this induced subgraph is fixed, pointwise, under each automorphism of $G$ preserving the labeling.  Thus each automorphism  of $G$ preserving the labeling, fixes the vertex $v_0$, each vertex of $N_G(v_0)$,  and so preserves the distance of $v_0$. We show that all vertices at distance $i$, $i\geq 1$, from  $v_0$ are fixed    under each automorphism of $G$ preserving the labeling, by induction on $i$. 

Let $x$ be  vertex at distance $i$, $i\geq 1$, of $v_0$ which is fixed  under each automorphism of $G$ preserving the labeling. Since the label of vertices adjacent to $x$ in $G$ which are at distance $i+1$ from $v_0$ are distinct, so all of them are fixed, too. This argument can be used for all vertices at distance $i$ from $v_0$, and so we conclude that all vertices  at distance $i+1$ from $v_0$  are fixed. Therefore the identity automorphism is the only automorphism of $G$ preserving the labeling and so $D(G)\leq \Delta -1$. \qed

\begin{figure}
	\begin{center}
		\includegraphics[width=0.3\textwidth]{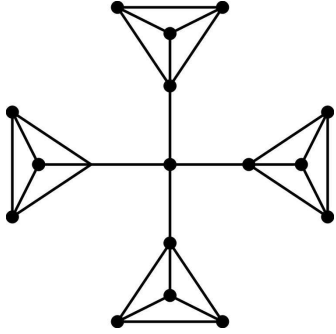}
		\caption{\label{exception} A graph with $\Delta (G)=\chi (G)= D(G)= 4$.}
	\end{center}
\end{figure}

The condition $\Delta (G) \geq 5$ is necessary in Theorem \ref{ubboundplana}. For instance, see Figure \ref{exception} for a graph with $\Delta (G)=\chi (G)= D(G)= 4$. 
It is known that a connected planar graph with chromatic number $\chi (G)=4$ has clique number 4, so  the following result is an immediate  consequence of Theorem \ref{ubboundplana}.

\begin{corollary}
	Let $G$ be a connected planar graph with chromatic number $\chi (G)=4$ and maximum degree $\Delta \geq 5$. Then $D(G)\leq \Delta -1$.
\end{corollary}

\subsection{Maximal outerplanar graphs}

\begin{figure}
	\begin{center}
		\includegraphics[width=0.8\textwidth]{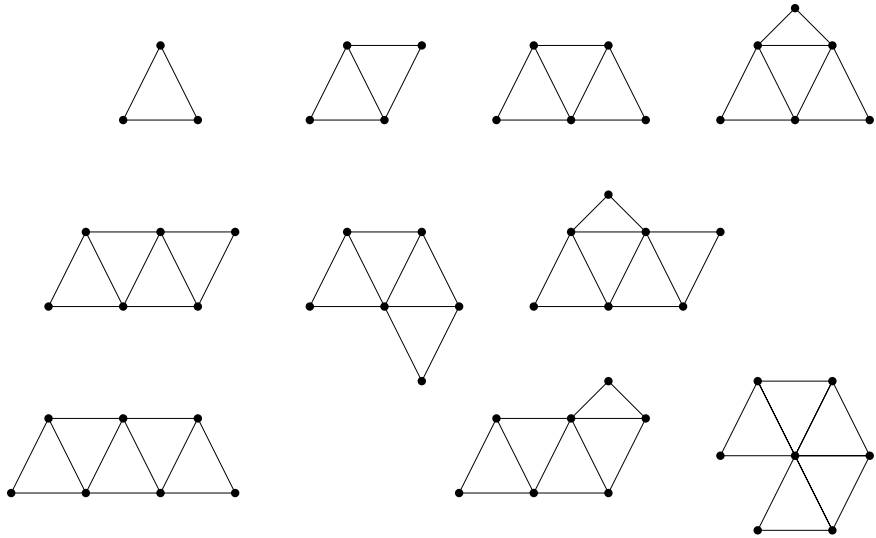}
		\caption{\label{mop}Maximal outerplanar graph with seven or less vertices.}
	\end{center}
\end{figure}

A maximal outerplanar graph, or briefly a MOP graph, is a graph that is isomorphic to
a triangularization of a polygon. All MOPs can be constructed according to
the following recursive rule (see, for instance, \cite{Proskurowski}): (i) the triangle, $K_3$, is
a MOP, and (ii) a MOP with $n + 1$ vertices can be obtained from some MOP
$M$ with $n$ vertices $( n \geq 3 )$, by adding a new vertex adjacent to two
consecutive vertices on the Hamiltonian cycle of $M$. All MOPs with up to seven vertices are shown in Figure \ref{mop}.

\begin{theorem}\label{onejamilton}
	If $G$ is a connected graph of order $n\geq 6$ with exactly one Hamiltonian cycle, then $D(G)\leq 2$.
\end{theorem}
\proof Let $V(G)=\{v_1,\ldots , v_n\}$. Since there is only one cycle, say $C$, with consecutive vertices $v_1,\ldots , v_n$, so the restriction of any automorphism of $G$ to $C$ is an automorphism of $C$.  We label the vertices of $G$ with two labels 1 and 2 such that  the cycle $C$ with consecutive vertices $v_1,\ldots , v_n$ has been distinguishingly. This labeling is a distinguishing labeling for $G$, Since if $f$ is an automorphism of $G$ preserving this labeling, then we can consider $f$ as an automorphism of $C$ preserving its labeling, and so $f$ is the identity automorphism, since we labeled $C$ distinguishingly. \qed

\begin{corollary}\label{numbmop}
	If $G$ is a maximal outerplanar graph, then $D(G)\leq 2$, except for $K_3$.
\end{corollary}
\proof It is known  that every maximal outerplanar graph  has exactly one Hamiltonian cycle. By Figure \ref{mop}, it can be seen that the distinguishing number of maximal outerplanar graphs  with seven or less vertices is at most two, except for  $K_3$.  Hence the result is a direct consequence of Theorem \ref{onejamilton}.\qed

\begin{theorem}\label{indxmop}
	If $G$ is a maximal outerplanar graph, then  $D'(G)\leq 2$, except for  $K_3$. 
\end{theorem}
\proof By Figure \ref{mop}, it can be seen that the distinguishing index of maximal outerplanar graphs  with seven or less vertices is at most two, except for  $K_3$.  Hence the result is a direct consequence of Theorem \ref{thm10}.\qed

Corollary \ref{numbmop} and Theorem \ref{indxmop} are true  for $2$-connected outerplanar graphs, too, because these graphs have  exactly one Hamiltonian cycle.


\section{Halin graphs}
A  Halin graph is constructed by connecting the leaves of a tree into a cycle. The tree must have at least four vertices, none of which has exactly two neighbors; it should be drawn in the plane so that none of its edges cross, and the cycle connects the leaves in their clockwise ordering in this embedding. Every Halin graph is a Hamiltonian graph. We start with the following lemma.

\begin{lemma}\label{lemhalin}
	Let $T$ be  a tree of order at least four  with no vertex of degree two, and $G$ be its Halin graph. If $f$ is an automorphism of $G$ fixing the vertices  
	of the cycle  induced on the set of the  leaves of $T$, say $C$, then $f$ is the identity automorphism of $G$.
\end{lemma}
\proof
Let $f$ fixes the vertices of $C$, but moves  a vertex $x\in V(T)\setminus V(C)$ to $y$. Since $f$ fixes the vertices of $C$, so we can consider $f$ as an automorphism of $T$. On the other hand, since degree of $x$ in $T$ is at least 3, so the set of the nearest vertices of degree 1 to $x$ in $T$ is different from the set of the nearest vertices of degree 1 to $y$ in $T$. Now since $f$ maps $x$ to $y$, so $f$ moves some vertices of $C$ which is a contradiction.\qed

\begin{theorem}\label{disnumhalin}
	If  $T$ is a tree of order at least four  with no vertex of degree two, and $G$ be its Halin graph, then $D(G)\leq 4$. The equality holds for $K_{1,3}$.
\end{theorem}
\proof We label the vertices of the cycle induced on the set of the  leaves of $T$, say $C$, distinguishingly with three labels 1,2, and 3. Next we label the remaining unlabeled vertices of $G$ with a new label 4. This labeling is distinguishing. In fact, if $f$ is an automorphism of $G$ preserving the labeling, then the restriction of $f$ to vertices of $C$ is an automorphism of $C$. Since the cycle $C$ has been labeled distinguishingly, so this restriction is the identity, and thus $f$ is the identity automorphism of $G$, by Lemma \ref{lemhalin}.\qed

\begin{theorem}
	Let $T$ be a tree of order at least four   with no vertex of degree two, and $G$ be its Halin graph. If $T$ has no vertex of degree $3$, then $D(G)\leq 3$. The equality holds for wheels $W_4$ and $W_5$.
\end{theorem}
\proof Since $T$ has no vertex of degree 3, so the only vertices of degree 3 in $G$ is the vertices of the cycle inuced on the set of leaves of $T$, say $C$. Hence the automorphisms of $G$ map the vertices of $C$ to itself, setwise. Thus the automorphism group of $G$ is a subgroup of the automorphism group of $T$.  Now, we label the vertices of the cycle induced on the set of the  leaves of $T$, say $C$, distinguishingly with three labels 1,2, and 3. Next we label the remaining unlabeled vertices of $G$ with an arbitrary label, say 1. Then this labeling is distinguishing. In fact, if $f$ is an automorphism of $G$ preserving the labeling, then since the restriction of $f$ to vertices of $C$ is an automorphism of $C$ and the cycle $C$ has been labeled distinguishingly, so this restriction is the identity, and thus $f$ is the identity automorphism of $G$, by Lemma \ref{lemhalin}.\qed

\begin{theorem}
	Let $T$ be a tree of order at least four  with no vertex of degree two, and $G$ be its Halin graph. If $T$ has no vertex of degree $3$ and the number of leaves of $T$ is at least $6$, then $D(G)\leq 2$. 
\end{theorem}
\proof It is known that $D(C_n)=2$ for any $n\geq 6$. Thus if we label the vertices of the cycle induced on the set of the  leaves of $T$, say $C$, distinguishingly with two labels 1 and 2, and next we label the remaining unlabeled vertices of $G$ with an arbitrary label, say 1, then this labeling is distinguishing, by a similar argument as proof of Theorem \ref{disnumhalin}.\qed

We recall that a \textit{ traceable graph} is a graph that possesses a Hamiltonian path.
\begin{theorem}{\rm \cite{nord}}\label{thm10}
	If $G$ is a traceable graph of order $n \geq 7$, then $D'(G) \leq  2$.
\end{theorem}
The assumption $n \geq 7$ is substantial in this  theorem, because for example  $D'(K_{3,3}) = 3$. 
\begin{figure}
	\begin{center}
		\includegraphics[width=0.8\textwidth]{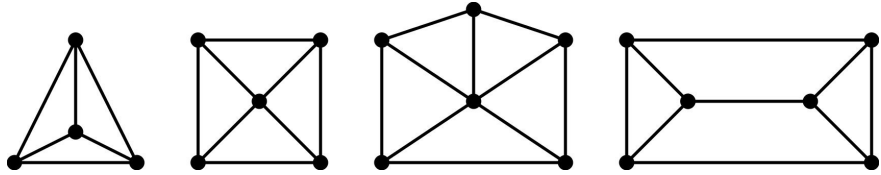}
		\caption{\label{hlinleq6}All Halin graphs of with six or less vertices.}
	\end{center}
\end{figure}

\begin{theorem}
	The distinguishing index of Halin graphs  is at most two, except for  $K_4$. 
\end{theorem}
\proof It can be seen that the distinguishing index of Halin graphs of order at most 6 is 2 except for $K_4$, see Figure \ref{hlinleq6}. For the Halin graphs of order at least 7, the result follows from Theorem \ref{thm10}.\qed

\section{ Mycielski graphs}
    Let the $n$ vertices of the given graph $G$ be $v_0$, $v_1,...,v_{n-1}$. The Mycielski graph $\mu(G)$ of $G$ contains $G$ itself as an isomorphic subgraph, together with $n+1$ additional vertices: a vertex $u_i$ corresponding to each vertex $v_i$ of $G$, and another vertex $w$. Each vertex $u_i$ is connected by an edge to $w$, so that these vertices form a subgraph in the form of a star $K_{1,n}$. In addition, for each edge $v_iv_j$ of $G$, the Mycielski graph includes two edges, $u_iv_j$ and $v_iu_j$. By above notations of vertices of  Mycielski graph $\mu(G)$ we have the following result. We recall that graphs with no pairs of vertices with the same open neighborhoods are called $R$-thin.
\begin{proposition}\label{propmyciel}
Let $G$ be an $R$-thin graph. If $f$ is an automorphism of $\mu(G)$ such that $f(w)=w$, then the restriction of $f$ to  each of sets $\{u_1, \ldots , u_n\}$ and $\{v_1,\ldots,v_n\}$ is isomorphic to an automorphism of $G$.
\end{proposition}
\proof If  $f(w)=w$, then since $f$ preserves adjacency relation, so $f$ maps the sets $\{u_1, \ldots , u_n\}$ and $\{v_1, \ldots , v_n\}$ to itself, setwise. By the same reasoning  the restriction of $f$ to $\{v_1, \ldots , v_n\}$ is an automorphism of $G$. In sequel, we want to show that the restriction of $f$ to the vertices  $\{u_1, \ldots , u_n\}$  is isomorphic to  the restriction of $f$ to $\{v_1, \ldots , v_n\}$. For this purpose we suppose that  $f(u_i)=u_{\sigma (i)}$ and $f(v_i)=v_{\tau (i)}$ where $\tau$ and $\sigma$ are   permutations of $\{1, \ldots , n\}$ and $1\leq i \leq n$, and show that   $\tau = \sigma$. Let $v_k$ be an arbitrary vertex of $G$ and $N_G(v_k)=\{v_{k_1}, \ldots , v_{k_t}\}$. Thus $N_{\mu (G)}(u_k)\setminus \{w\}=\{v_{k_1}, \ldots , v_{k_t}\}$, and hence $N_{\mu (G)}(u_{\sigma (k)})\setminus \{w\}=\{v_{\tau (k_1)}, \ldots , v_{\tau (k_t)}\}$. Then $N_{G}(v_{\sigma (k)}) =\{v_{\tau (k_1)}, \ldots , v_{\tau (k_t)}\}$. On the other hand  $N_{G}(v_{\tau(k)}) =\{v_{\tau (k_1)}, \ldots , v_{\tau (k_t)}\}$. Therefore $v_{\tau(k)}$ and $v_{\sigma(k)}$ have the same open neighbor in $G$. Since $G$ is an $R$-thin graph, so $\tau (k)= \sigma (k)$. Since $k$ is arbitrary, so we have the result.\qed

\begin{theorem}\label{distnumbmyc}
If  $G$ is   an $R$-thin graph of order $n\geq 2$, then $D(\mu (G))\leq D(G)+1$.
\end{theorem}
\proof We present a distinguishing vertex labeling of  $\mu (G)$ with $D(G)+1$ labels. Let $c:V(G)\rightarrow \{1,\ldots , D(G)\}$ be a distinguishing labeling of $G$. We extend $c$ to a $(D(G)+1)$-labeling $\overline{c}$ of $\mu (G)$ with $\overline{c}(v_i)=c(v_i)$, $\overline{c}(u_i)=c(v_i)$ for $1\leq i \leq D(G)$, and $\overline{c}(w)=0$. We claim that $\overline{c}$ is a distinguishing labeling of $\mu (G)$. In fact if $f$ is an automorphism of $\mu (G)$ preserving the labeling $\overline{c}$, then $f(w)=w$, and so the restriction of $f$ to vertices $v_1, \ldots , v_n$ is an automorphism of $G$ preserving the labeling $c$. Thus the restriction of $f$ to the vertices of $G$ is the identity, and so $f$ is the identity automorphism of $\mu (G)$ by Proposition \ref{propmyciel}.\qed

\begin{theorem}\label{distindmyc}
If  $G$ is   an $R$-thin graph of order $n\geq 3$ with no connected component $K_2$, then $D'(\mu (G))\leq D'(G)+1$.
\end{theorem}
\proof Let $c':E(G)\rightarrow \{1,\ldots , D'(G)\}$ be a distinguishing edge labeling of $G$. We define a $(D'(G)+1)$-labeling $\overline{c'}$ of $\mu (G)$ with $\overline{c'}(wu_i)=0$   for $1\leq i \leq D(G)$, $\overline{c'}(u_iv_j)=c'(v_iv_j)$ for $i,j \in \{1, \ldots , n\}$ for which $v_iv_j\in E(G)$,  and $\overline{c'}(e)=c'(e)$ for all $e\in E(G)$. The mapping $\overline{c'}$ is a distinguishing labeling of $\mu (G)$. In fact if $f$ is an automorphism of $\mu (G)$ preserving the labeling $\overline{c'}$, then $f(w)=w$, and so the restriction of $f$ to vertices $v_1, \ldots , v_n$ is an automorphism of $G$ preserving the labeling $c'$. Thus the restriction of $f$ to the vertices of $G$ is the identity, and so $f$ is the identity automorphism of $\mu (G)$ by Proposition \ref{propmyciel}. \qed

 Applying the Mycielskian repeatedly, starting with a graph with a single edge, produces a sequence of graphs $M_i = \mu(M_{i-1})$, also sometimes called the Mycielski graphs. The first few graphs in this sequence are the graph $M_2 = K_2$ with two vertices connected by an edge, the cycle graph $M_3 = C_5$, and the Gr\"otzsch graph with 11 vertices and 20 edges.
It is clear that $M_i$ is an $R$-thin graph if and only if $M_{i-1}$ is an $R$-thin graph. Since $M_2$ is an $R$-thin graph, so $M_i$'s are $R$-thin graphs for any $i\geq 2$.

 \begin{figure}
 	\begin{center}
 		\includegraphics[width=0.8\textwidth]{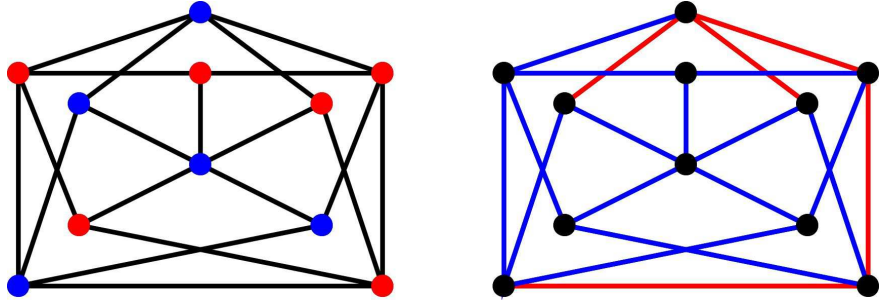}
 		\caption{\label{coloringM4}A 2-distinguishing labeling of vertices and edges of $M_4$.}
 	\end{center}
 \end{figure}

\begin{theorem}\label{induct}
For any $i\geq 5$, $D(M_i)\leq D(M_{i-1})$ and  $D'(M_i)\leq D'(M_{i-1})$. 
\end{theorem}
\proof The degree of vertex $w$ in $M_i$ is $|V(M_{i-1})|$, and $w$ is the only vertex with this  degree, so $w$ is fixed under each automorphism of $M_i$. We first show that $D(M_i)\leq D(M_{i-1})$. We label the vertices of isomorphic subgraph $M_{i-1}$ in $M_i$, i.e., $v$'s, distinguishingly with $D(M_{i-1})$ labels, and label the corresponding vertices to the vertices $M_{i-1}$ in $M_i$, i.e., $u$'s, exactly the same as $v's$ and label $w$ with an arbitrary label. Then by the same reasoning as Theorem \ref{distnumbmyc} we conclude that this labeling is distinguishing. Therefore  $D(M_i)\leq D(M_{i-1})$.   

For the second part of Theorem, we label the incident edges to $w$ with an arbitrary label, say 1, and label the remaining edges of $M_i$ exactly the same as Theorem \ref{distindmyc}. Then by the same reasoning  we conclude that this labeling is distinguishing. Therefore  $D'(M_i)\leq D'(M_{i-1})$. \qed

\begin{corollary}
For $i\geq 4$, $D(M_i)=D'(M_i)=2$. In particular, $D(M_3)=D'(M_3)=3$ and $D(M_2)=2$.
\end{corollary}
\proof It is easy to see that $D(M_2)=2$ and $D(M_3)=D'(M_3)=3$. In Figure \ref{coloringM4} we presented a 2-distinguishing labeling of vertices and edges of $M_4$. Now the result follows from Theorem \ref{induct} and the  induction on $i$. \qed

We end this paper by the following conjecture:
\begin{conjecture}
 Let $G$ be  a connected graph of order $n\geq 3$. Then $D(\mu (G))\leq D(G)$ and $D'(\mu (G))\leq D'(G)$, except for a finite number of graphs.
\end{conjecture}

\end{document}